\font \svntnrm =cmr17
\font \svntnmi =cmmi12 scaled \magstep 2     
\font \frtnrm =cmr12 scaled \magstep 1       
\font \frtnmi =cmmi12 scaled \magstep 1     
\font \elvnrm =cmr10 scaled \magstephalf
\font \elvnmi =cmmi10 scaled \magstephalf
\def \titlefont {%
        \textfont0 =\svntnrm \scriptfont0 =\frtnrm
                \scriptscriptfont0 =\elvnrm
        \textfont1 =\svntnmi \scriptfont1 =\frtnmi
                \scriptscriptfont1 =\elvnmi
        \def \rm  {\fam0 \svntnrm}%
        \def \mit {\fam1 }%
        \rm }%

\font\sc=cmcsc10
\def \title #1{ {\titlefont \centerline {#1}}}
\def \address #1 { {\sl \baselineskip =10pt \vskip 2truemm
                      {\halign {\centerline {##}\cr 
                          #1}}\vskip 2truemm}}

\def \and {\&\ }


%
%
\def \section #1#2{ \goodbreak 
\bigskip \centerline {\bf #1{\def \arg
    {#1}\ifx \arg \empty \else .\fi}\quad #2}\nobreak \bigskip \xdef
    \thissection {#1}}

\magnification =\magstep 1
\overfullrule =0pt
%
%
\newwrite \reffile 
\immediate \openout \reffile =\jobname.ref
\newcount \refnum \refnum =0
\def \ref #1: #2{\global \advance \refnum by 1 \expandafter \xdef 
    \csname ref#1\endcsname {\the \refnum }\immediate \write \reffile 
    {\noexpand \item {[\csname ref#1\endcsname ]}{#2}}}
%
%
\ref Concrete: {R.L. Graham, D.E. Knuth, O. Patashnik,
 {\sl Concrete Mathematics}, Addison-Wesley, 1994.}
\ref AdvCom: {L. Comtet, {\sl Analyse Combinatoire},2, Presses Univ.
de France, 1970.}
\ref deB: {N. G. de Bruijn, {\sl Asymptotic Methods in Analysis\/},
North-Holland, 1961.}
\ref Cara: {C. Carath\'eodory, {\sl Theory of Functions of a Complex
Variable}, Chelsea, 1954.}
\ref Comtet: {L. Comtet, {\sl C. R. Acad. Sc. Paris\/}, 270, 1970, 
p. 1085--1088.}
\ref Wpap: {R. M. Corless, G. H. Gonnet, D. E. G. Hare, 
D. J. Jeffrey and D. E. Knuth, ``The Lambert $W$ Function'',
Advances in Computational Mathematics, to appear.}
\ref David: {F.N. David and D.E. Barton, {\sl Combinatorial Chance},
Hafner, 1962.}
\ref Niel: {N. Nielsen, {\sl Handbuch der Theorie der Gammafunktion},
Teubner, 1906.}
\ref Rior: {J. Riordan, {\sl An Introduction to Combinatorial
Analysis}, Wiley, 1958.}
\immediate \closeout \reffile 
%
%



\bigskip 
\title {On the inversion of $y^\alpha e^y$ in terms of associated Stirling numbers}
\medskip 
\centerline { by}
\smallskip 
\centerline {D.~J.~Jeffrey${}^1$, R.~M.~Corless${}^1$, 
D.~E.~G.~Hare${}^2$, and D.~E.~Knuth${}^3$}
\bigskip 
{\narrower
The function $y=\Phi_\alpha(x)$, the solution of $y^\alpha e^y=x$
for $x$ and $y$ large enough,
has a series expansion in terms of $\ln x$ and $\ln\ln x$,
with coefficients given in terms of Stirling cycle numbers.
It is shown that this expansion converges for $x>(\alpha e)^\alpha$
for $\alpha \ge 1$.
It is also shown that new expansions 
can be obtained for $\Phi_\alpha$ in terms of associated Stirling
numbers. The new expansions converge more rapidly
and on a larger domain.   \smallskip}

\noindent 1. {\sc Stirling numbers} ---
Stirling cycle numbers ${ n\brack m}$ are defined [\refConcrete] by
$$ \ln^m(1+z) = m! \sum_n (-1)^{n+m} { n\brack m} {z^n\over n!}
  \ .\leqno(1a) $$
The numbers $(-1)^{n+m}{n\brack m}$ are also called Stirling numbers of the first kind [\refNiel].
Stirling subset numbers ${n\brace m}$, also called Stirling numbers
of the second kind, are defined by
$$ \left( e^z-1\right)^m = m! \sum_n {n\brace m} {z^n\over n!}
  \ ,\leqno(1b) $$
and 2-associated Stirling subset numbers ${n\brace m}_{\ge 2}$ are
defined by 
[\refAdvCom, exercise 5.7; \refDavid, p. 296; \refRior, \S 4.5]
$$ \left( e^z-1-z\right)^m = m!\sum_n {n\brace m}_{\!\ge 2} {z^n\over n!}
  \ .\leqno(1c) $$

\noindent 2. {\sc Solution by Comtet of} $y^\alpha e^y =x$. ---
The fixed real $\alpha$,
we let $\Phi_\alpha (x) $ be the value of $y$ that is the unique 
positive solution of the equation $y^\alpha e^y = x$.
If $\alpha$ is negative, then $ y > -\alpha $ and 
$x> e^{-\alpha} (-\alpha)^\alpha$.
An asymptotic expansion for $\Phi_\alpha(x)$,
in terms of Stirling cycle numbers and 
the quantities $L_1=\ln x$ and $L_2=\ln\ln x$,
is given in the following theorem [\refdeB,\refComtet].

Theorem 1. --- {\it With the preceding notation, 
the function $\Phi_\alpha(x)$ has the following series development, 
convergent if $x$ is large enough.}
$$ \Phi_\alpha(x) = L_1 - \alpha L_2 +\alpha
  \sum_{n\ge 1} {\alpha^n\over L_1^n}
      \sum_{m=1}^n (-1)^{n+m} {n\brack n-m+1 } {L_2^m\over m!}
 \ .\leqno(2a)$$

{\it Proof.} -- We recall some details of the proof given in
[\refComtet] for use below.
We introduce a function $w(x)$ defined by
$$ y=\Phi_\alpha(x) = L_1-\alpha L_2 +\alpha w \ ,\leqno(2b) $$
which satisfies 
$$ 1-e^{-w} +\sigma w-\tau=0,\qquad \sigma = {\alpha\over L_1},\qquad
\tau = \alpha {L_2\over L_1}
=\sigma\ln\left({\alpha\over\sigma}\right)
\ . \leqno(2c) $$
By the Lagrange Inversion Theorem [\refCara], $w$ has the expansion
$$ w= \sum_{m\ge 1} {\tau^m\over m!} 
      \sum_{l\ge 0} (-1)^l {l+m \brack l+1} \sigma^l \>.
        \leqno(2d) $$
One converts from $\sigma$ and $\tau$ back to $L_1$ and $L_2$ to
complete the theorem.

Since the domain of convergence of $(2a)$ is described only 
as `$x$ large enough'
by de~Bruijn and Comtet, we give a stronger statement in the next theorem.

Theorem 2. --- 
{\it For $\alpha \ge 1$, the series (2a) is convergent for 
$x>(\alpha e)^\alpha$, while for $\alpha <1$ it is convergent for
$x>e$.}

{\it Abbreviated Proof.} -- We let $f(w)=\sigma w -\tau$ 
and $g(w)=1-e^{-w}$. 
For $\alpha>1$ and $x>(\alpha e)^\alpha$, 
we define $\delta>0$ by $\delta=1-\ln (\alpha e)^\alpha /\ln x$, and
then $\sigma=(1-\delta)/(1+\ln\alpha)$. We also set
$w_0=\ln(1+\ln\alpha)$.
On the contour consisting of the lines $\Re(w)=w_0+\delta$, 
$\Im(w)=\pm 2\delta^{1/2}$, and $\Re(w) = -2$,
one can show that $\vert g\vert >\vert f\vert$,
and therefore $f+g$
has only one root within the contour, by Rouch\'e's theorem.
Using Cauchy's theorem to express this root as an integral around the
contour, we establish the convergence of $(2a)$  
by expanding the integrand as a series in $f/g$ and integrating
term by term [\refdeB].
For $\alpha<1$, the contour must remain the same as that for
$\alpha=1$.
\smallskip
\noindent 3. {\sc A New Expansion.} ---
In view of the relation
$$ \Phi_\alpha(x)
=\alpha \Phi_1\left( {x^{1/\alpha} \over \alpha}\right)
=\alpha W\left( {x^{1/\alpha} \over \alpha}\right)\ ,$$
where $W$ is the Lambert $W$ function
[\refWpap],
we shall simplify our equations by considering
only the case $\alpha=1$ from now on.
By changing to the variable $\zeta = 1/(1+\sigma)$,
we obtain a new series for $W=\Phi_1$
that converges on a wider domain than does $(2a)$.

Theorem 3. --- {\it  With the preceding notation, 
$W$ has the series development }
$$ W(x) = L_1 - L_2 +
  \sum_{m\ge 1} {\tau^m\over m!}\, \sum_{p=0}^{m-1}
(-1)^{p+m-1} \zeta^{p+m}{ p+m-1 \brace p}_{\!\ge 2}
 \>, \leqno(3a)   $$
{\it and this is convergent for} $x \ge 2$.

{\it Proof.} ---  
Into $(2c)$, we substitute $\sigma=1/\zeta -1$ and obtain
$$ \tau +e^{-w}-1+w - w/\zeta = 0 \ .\leqno(3b)$$
To invert this using the Lagrange Inversion Theorem, we introduce the
operator $[w^p]$ to represent the coefficient of $w^p$ in a series
expansion in $w$, and obtain
$$\eqalign{ 
 w &= \sum_{n\ge 1} {\zeta^n\over n} \left[ w^{n-1}\right]
     \left( \tau + e^{-w} -1+w\right)^n\ ,\cr
&= \sum_{n\ge 1} {\zeta^n\over n} \left[ w^{n-1}\right]
\sum_m {n \choose m} \tau^m \left( e^{-w}-1+w\right)^{n-m} \cr
&= \sum_{n\ge 1} (-1)^{n-1} \zeta^n \sum_m {\tau^m\over m!}
 { n-1\brace n-m}_{\!\ge 2}\cr } $$
which can be rearranged to obtain the theorem.

To prove convergence, we let $f(w)=\zeta(e^{-w}-1+w)+\tau\zeta$ and 
$g(w)= -w$.
On the rectangular contour bounded by the four lines
$\Re(w)=2 $, $\Im(w)=\pm 2$ and $\Re(w) = -1$,
it is simple to show that $\vert f\vert <\vert g\vert $ for all
$x\in[2,e]$. Hence the series converges there. Since $(3a)$ is equivalent
to $(2a)$ for $x>e$ because of the relation
$$ {l\brack m} = \sum_{p=0}^{l-m} (-1)^{p+l-m}
 {p+l-m\brace p}_{\!\ge 2}
{ p+l-1 \choose p+l-m} \ ,\leqno(3c)   $$
the theorem follows.


\smallskip
\noindent 4. {\sc Expansions using new variables.} ---
Two further series developments can be obtained by introducing the
variables $L_\tau= \ln(1-\tau)$ and $\eta=\sigma/(1-\tau)$.

Theorem 4.---{\it With the preceding notation, $W$ has the series development}
$$ W(x) = L_1- L_2 - L_\tau
   - \sum_{n\ge 1}(-\eta)^n
       \sum_{m=1}^n (-1)^{m+1} {n\brack n-m+1} {L_\tau^m\over m!}\ .
\leqno(4a)$$

{\it Proof.}---
We set $w=v-L_\tau$ in $(2c)$ and obtain, after rearranging,
$$ 1- e^{-v} + {\sigma\over 1-\tau} v = {\sigma\over 1-\tau} L_\tau
\ .\leqno(4b) $$
This equation has exactly the form of $(2c)$ itself, 
and therefore the expansion for $v$ can be obtained from $(2d)$
by replacing $\sigma$ with $\sigma/(1-\tau)$ and $\tau$ with
$\sigma L_\tau/(1-\tau)$. The theorem then follows by
rearrangement.

The expansion $(4a)$ converges more slowly than $(2a)$,
but when we transform it using the methods of theorem 3,
we obtain a very rapidly convergent expansion, 
as we show in section 5.

Theorem 5.---{\it With the above notation, 
$W(x)$ has the development}
$$ W(x) = L_1 - L_2 - L_\tau+
   \sum_{m\ge 1}{1\over m!} L_\tau^m \eta^m 
     \sum_{p=0}^{m-1}(-1)^{p+m-1}{ p+m-1 \brace p}_{\!\ge 2}
 { 1 \over (1+\eta)^{p+m} }\>,
\leqno(4c) $$

{\it Proof.} --- The proof follows exactly that of Theorem 3.

The process of generating series in new variables can be continued.
If $w(\sigma,\tau)$ satisfies $(2c)$, then Theorem 4 is equivalent to
the identity
$$ w(\sigma,\tau) = - \ln(1-\tau) + w\left({\sigma\over 1-\tau},{\sigma
\ln(1-\tau)\over 1-\tau}\right)\ , \leqno(4d) $$
which clearly can be applied repeatedly.
\smallskip
\noindent 5. {\sc Rate of convergence.} ---
We consider the accuracy obtained by truncating each of the
series $(2a)$, $(3a)$ and $(4c)$ at $N-1$ terms. 
Since the series are asymptotic series, the error terms for $x$ large are
respectively $O(L_2^N/L_1^N )$ for $(2a)$ and $(3a)$ and
$O(L_2^N/L_1^{2N})$ for (4c), so $(4c)$ is clearly better.
In addition to being asymptotic, however, the series are
absolutely convergent, and can be used for relatively small values
of $x$.
We observe that $\tau=L_\tau=0$ at $x=e$, 
and hence the infinite sums in $(2a)$, $(3a)$ and $(4c)$ are zero there.
Thus any truncated series will be exact at $x=e$ and
asymptotically correct as $x\to \infty$, implying that the error
will have a maximum at some $x>e$.

However, although $(2a)$ is correct at $x=e$, its derivative
does not converge there.
In contrast, $(3a)$ and $(4c)$ give finite sums at $x=e$ for all
derivatives.
To put it another way, taking $N$ terms of $(3a)$ or $(4c)$ and expanding
about $x=e$ gives $N$ terms of the Taylor series for $\Phi_1(x)$
about $x=e$.
Both $(3a)$ and $(4c)$ are much more accurate than $(2a)$ near this point.
Numerical experiments confirm these results.

We conjecture that (2a) and (4c) converge for all $x>1$.
%
%
\section {}{References}
{\frenchspacing
\input \jobname.ref
}
\smallskip 
\address {${}^1$ Department of Applied Mathematics \cr 
          The University of Western Ontario        \cr 
          London, CANADA, N6A 5B7              \cr 
          \cr
          ${}^2$ Symbolic Computation Group \cr
          University of Waterloo \cr
	  Waterloo, CANADA, N2L 3G1 \cr
          \cr
          ${}^3$ Department of Computer Science \cr
          Stanford University \cr
          Stanford, USA, 94305-2140 \cr
         }
\vfill\eject
\bye